\documentclass[12pt]{amsart}
\usepackage{amsmath,amssymb,amsfonts,amsthm}
\usepackage{fullpage}

\usepackage{color,graphicx}

\newtheorem{thm}{Theorem}[section]
\newtheorem{prop}[thm]{Proposition}
\newtheorem{lemma}[thm]{Lemma}

\theoremstyle{definition}
\newtheorem{defn}[thm]{Definition}
\theoremstyle{remark}
\newtheorem{rmk}[thm]{Remark}

\newtheorem{ex}[thm]{Example}

\newcommand{\abs}[1]{\left\lvert#1\right\rvert}

\newcommand{\Z}{\mathbb{Z}}

\usepackage[latin1]{inputenc}
\usepackage{subfigure}
\usepackage[numbers]{natbib}

\graphicspath{{./figures/}}

\author{Michael Chmutov}
\title{Type $A$ molecules are Kazhdan-Lusztig}
\address{Department of Mathematics, University of Michigan, 2074 East Hall, 530 Church St., Ann Arbor, MI 48109}

\begin{document}
\begin{abstract}
Let $(W, S)$ be a Coxeter system. A $W$-graph is an encoding of a representation of the corresponding Iwahori-Hecke algebra. Especially important examples include the $W$-graph corresponding to the action of the Iwahori-Hecke algebra on the Kazhdan-Lusztig basis, as well as this graph's strongly connected components (cells). In 2008, Stembridge identified some common features of the Kazhdan-Lusztig graphs and gave a combinatorial characterization of all $W$-graphs that have these features. He conjectured, and checked up to $n=9$, that all such $A_n$-cells are Kazhdan-Lusztig cells. The current paper provides a first step toward a potential proof of the conjecture. More concretely, we prove that the connected subgraphs of $A_n$-cells consisting of simple (i.e. directed both ways) edges are the same as in the Kazhdan-Lusztig cells.
\end{abstract}
\maketitle

\section{Introduction}

Let $(W, S)$ be a Coxeter system. A $W$-graph is a directed graph with additional structure that encodes a representation of the corresponding Iwahori-Hecke algebra.
In the paper \cite{KL}, such graphs were introduced for the regular representation, and it was shown that the strongly connected components (called ``cells'') also yield representations. 
Stembridge identified several common features of the Kazhdan-Lusztig graphs, namely, they are bipartite, (nearly) edge-symmetric, and their edge weights are non-negative integers (collectively he called these features ``admissibility''). He proceeded to describe, via four combinatorial rules, when an admissible graph is a $W$-graph (\cite{st_adm}). One hopes that the characterization will allow 
one to construct the Kazhdan-Lusztig cells without having to compute Kazhdan-Lusztig polynomials (a notoriously difficult task). A piece of evidence suggesting that the definition of a general admissible $W$-cell approximates a Kazhdan-Lusztig cell is a more recent result of Stembridge that there are only finitely many admissible $W$-cells for each finite $W$ (\cite{st_fin}).

There are no known examples of admissible $A_n$-cells besides the Kazhdan-Lusztig cells (Stembridge checked it up to $n=9$). A possible strategy of proof is as follows:
\begin{enumerate}
\item An $A_n$-cell is a strongly connected directed graph. It turns out that the \emph{simple}, i.e. directed both ways, edges are considerably easier to understand than the rest of the edges. Consider the induced subgraphs which are connected via simple edges (of course these are strongly connected on their own, but a cell may contain several of them). The subgraphs satisfy combinatorial rules slightly weaker than those satisfied by a cell; a graph satisfying these rules is called a \emph{molecule}. The first step is to show that the simple edges of any $A_n$-molecule appear in the Kazhdan-Lusztig graph.

\item It is known that a Kazhdan-Lusztig $A_n$-cell is connected via simple edges, and these edges are well understood (they are called \emph{dual Knuth moves}). The second step is to prove that no cell may have multiple molecules. The fact that no two Kazhdan-Lusztig $A_n$-molecules may be connected inside a cell has been checked for $n\leqslant 12$ (\cite{st_per}).

\item The last part is to prove that there can be only one $A_n$-graph with a given underlying molecule. For Kazhdan-Lusztig molecules this has been checked 
for $n\leqslant 13$ (\cite{st_per}).
\end{enumerate}

In this paper we complete the first part of the above program. Together with the above computations, this result implies that all $A_n$-cells up to $n=12$ are Kazhdan-Lusztig. The main ingredient of the proof is the axiomatization of graphs on tableaux generated by dual Knuth moves (\cite{assaf_degs}). Five of the axioms follow easily from the molecules axioms, but the last one presents a challenge. Recently Roberts suggested a revised version of the last axiom (\cite{roberts_degs}). Using it one can give a short computerized proof of our result. 

The paper is structured as follows. Section \ref{sec:mols} introduces the $W$-molecule world. Section \ref{sec:degs} discusses Assaf's dual equivalence graphs and connects them to molecules. The last section is devoted to the proof of the main theorem that the simple part of any $A_n$-molecule is a dual equivalence graph. 

\section{Molecules and arc transport}
\label{sec:mols}

This section summarizes the required $W$-molecules terminology as described in \cite{st_adm,st_atlas_mol}. 

Let $(W,S)$ be a finite simply-laced Coxeter system.  
The papers are mostly concerned with $W$-graphs, i.e. graphs that encode certain representations of the corresponding Iwahori-Hecke algebra. 
It turns out that the simple (i.e. directed both ways) edges of these graphs are much easier to understand than other edges. Thus we consider induced 
subgraphs connected by simple edges. These subgraphs are not necessarily $W$-graphs (i.e. they do not encode representations), but they satisfy certain 
combinatorial rules which are slightly weaker than Stembridge's $W$-graph rules. More specifically, to get the $W$-graph rules one needs to omit the conditions 
on $\tau$-invariants of $u$ and $v$ in the (LPR2) and (LPR3) sections of Definition \ref{defn:mols} below.
We begin this paper by formalizing the definitions and presenting the 
rules. See the above papers for more details regarding the relationship of $W$-molecules with $W$-graphs and representation theory. 
A significant part of this section extends to multiply-laced types; again, see the above two papers.

\subsection{Definitions and Axioms}
An \emph{admissible $S$-labeled graph} is a tuple $G=(V,m,\tau)$, where $V$ is a set (vertices), $m:V\times V\to\Z^{\geqslant 0}$, and $\tau:V\to 2^S$ such that
\begin{enumerate}
\item as a directed graph (with edges given by pairs of vertices with non-zero $m$ value), $G$ is bipartite,
\item if $\tau(u)\subseteq\tau(v)$ then $m(u, v)=0$,
\item if $\tau(u)$ and $\tau(v)$ are incomparable, then $m(u,v) = m(v,u)$.
\end{enumerate} 
The function $\tau$ will be referred to as the \emph{$\tau$-invariant}. We will most of the time omit the word ``admissible'' since we consider no other $S$-labeled graphs. 
By a \emph{simple edge} we mean a pair of vertices $(v_1, v_2)$ such that neither $m(v_1,v_2)$ nor $m(v_2,v_1)$ are 0. 
In diagrams we draw these as undirected edges (see the left side of Figure \ref{fig:a4molsdegs}).
By an \emph{arc} $v_1\to v_2$ we mean a pair of vertices $(v_1, v_2)$ such that $m(v_1,v_2)\neq 0$, but $m(v_2,v_1)=0$. Notice that if 
$u\to v$ is an arc, then $\tau(u)\supset\tau(v)$. If $(u,v)$ is a simple edge then $\tau(u)$ and $\tau(v)$ are incomparable, 
and $m(u,v) = m(v,u)$.

We refer to an edge of the Coxeter graph of $(W,S)$ as a \emph{bond} to distinguish it from edges of $S$-labeled graphs. A simple edge $(u, v)$ \emph{activates} a bond $(i,j)$ if precisely one of $\tau(u)$ and $\tau(v)$ contains $i$, and precisely the other one contains $j$.

For distinct $i,j\in S$, a directed path (possibly involving simple edges) $u\to v_1\to v_2\to\dots\to v_{r-1}\to v$ in $G$ is \emph{alternating of type $(i,j)$} if
\begin{itemize}
\item $i,j\in\tau(u)$ and $i,j\notin\tau(v)$,
\item $i\in\tau(v_k), j\notin\tau(v_k)$ for odd $k$,
\item $i\notin\tau(v_k), j\in\tau(v_k)$ for even $k$.
\end{itemize} 
Let 
$$N_{ij}^r(G;u,v):=\sum_{v_1,\dots,v_{r-1}} m(u,v_1)m(v_1,v_2)\dots m(v_{r-1},v),$$
where the sum is over the set of alternating paths of type $(i,j)$ from $u$ to $v$.

\begin{defn}
\label{defn:mols}
An $S$-labeled graph is called a \emph{molecular graph} if it satisfies
\begin{enumerate}
\item[(SR)] If $(u,v)$ is a simple edge then $m(u,v) = m(v,u) = 1$. Thus we will omit the weights of simple edge in our diagrams. 
\item[(CR)] If $u\to v$ is an edge, i.e. $m(u,v)\neq 0$, then every $i\in\tau(u)\setminus\tau(v)$ is bonded to every $j\in\tau(v)\setminus\tau(u)$.
\item[(BR)] Suppose $(i, j)$ is a bond in the Coxeter graph of $(W,S)$. Any vertex $u$ with $i\in\tau(u)$ and $j\notin\tau(u)$ is adjacent to precisely one edge which activates $(i, j)$.
\item[(LPR2)] For any $i,j\in S$ for any $u,v\in V$ with $i,j\in\tau(u)$, $i,j\notin\tau(v)$ and $\tau(v)\setminus\tau(u)\neq\varnothing$, we have
$$N_{ij}^2(G;u,v)=N_{ji}^2(G;u,v).$$
\item[(LPR3)] Let $k,i,j,l\in S$ be a copy of $A_4$ in the Coxeter graph: $k-i-j-l$. For any $u,v\in V$ with $i,j\in\tau(u)$, $i,j\notin\tau(v)$, $k,l\notin\tau(u)$, $k,l\in\tau(v)$, we have
$$N_{ij}^3(G;u,v)=N_{ji}^3(G;u,v).$$
\end{enumerate}
\end{defn}
The rules are called, respectively, simplicity rule, compatibility rule, bonding rule, and local polygon rules.

\begin{defn}
An $S$-labeled graph is called a \emph{molecule} if it is a molecular graph, and there is a path of simple edges between any pair of vertices.
\end{defn}

\begin{ex}
\label{ex:a3mols}
It is easy to classify all the $A_3$-molecules. Because of admissibility, a vertex whose $\tau$-invariant is $\varnothing$ cannot be connected to any other vertex by a simple edge.
Similarly for a vertex whose $\tau$-invariant is $\{1,2,3\}$. 

Suppose we have a vertex $v_1$, whose $\tau$-invariant is $\{1\}$. By BR, it is connected by a simple edge to a vertex $v_2$ whose $\tau$-invariant contains $2$, but not $1$. 
By CR, $3\notin\tau(v_2)$, and hence $\tau(v_2) =\{2\}$. By BR, $v_2$ is connected by a simple edge to a vertex $v_3$ whose $\tau$-invariant contains $3$, but not $2$. 
We already know $v_3 \neq v_1$. By BR, $\tau(v_3) = \{3\}$. There are no other simple edges possible, and this is a complete molecule. The same analysis works for $v_1$ having
$\tau$-invariants of $\{3\},\{1,2\}, \{2,3\}$.

Suppose we have a vertex $v_1$, whose $\tau$-invariant is $\{2\}$. By BR, it is connected by a simple edge to a vertex $v_2$ whose $\tau$-invariant contains $1$, but not $2$.
The case of $\tau(v_2)=\{1\}$ was described above, so the only choice is $\tau(v_2) = \{1,3\}$. This yields a complete molecule. The same argument works for  $v_1$ having
$\tau$-invariant of $\{1,3\}$.

This completes the classification:
$$\includegraphics[width=6cm]{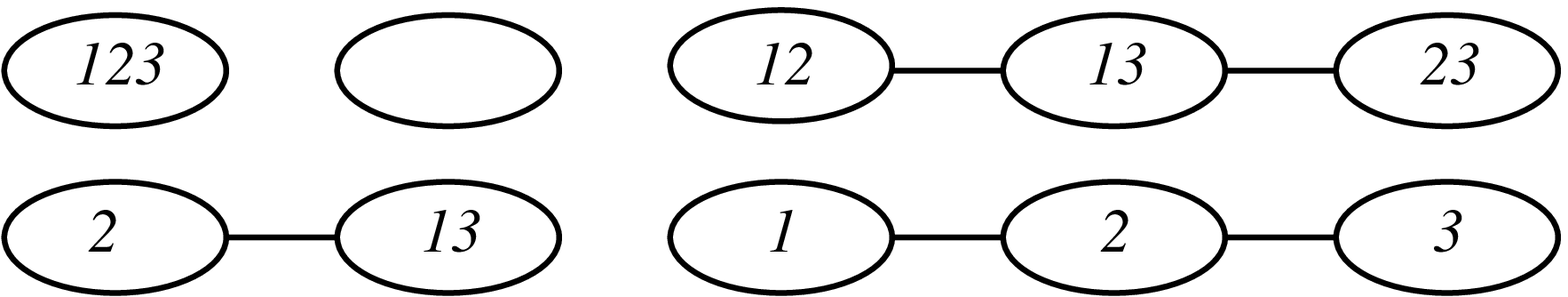}$$
\end{ex}

\begin{ex}
\label{ex:a4mols}
The classification of $A_4$-molecules is slightly more involved (see \cite{st_adm}). The result is shown on the left side 
of Figure \ref{fig:a4molsdegs}.

\begin{rmk}
\label{rmk:lprpaths}
We would like to comment on the structure of alternating paths involved in the local polygon rules. A priori only the first and the last edges of an alternating path could be arcs. In fact, the additional assumptions on the $\tau$-invariants of the starting and ending vertices force one of these edges to be simple. So any alternating path involved in the local polygon rules contains at most one arc. 
\end{rmk}

\begin{figure}
\caption{Molecules and dual equivalence graphs for type $A_4$}
\label{fig:a4molsdegs}
\bigskip
\begin{tabular}{l|r}
\includegraphics[height=8cm]{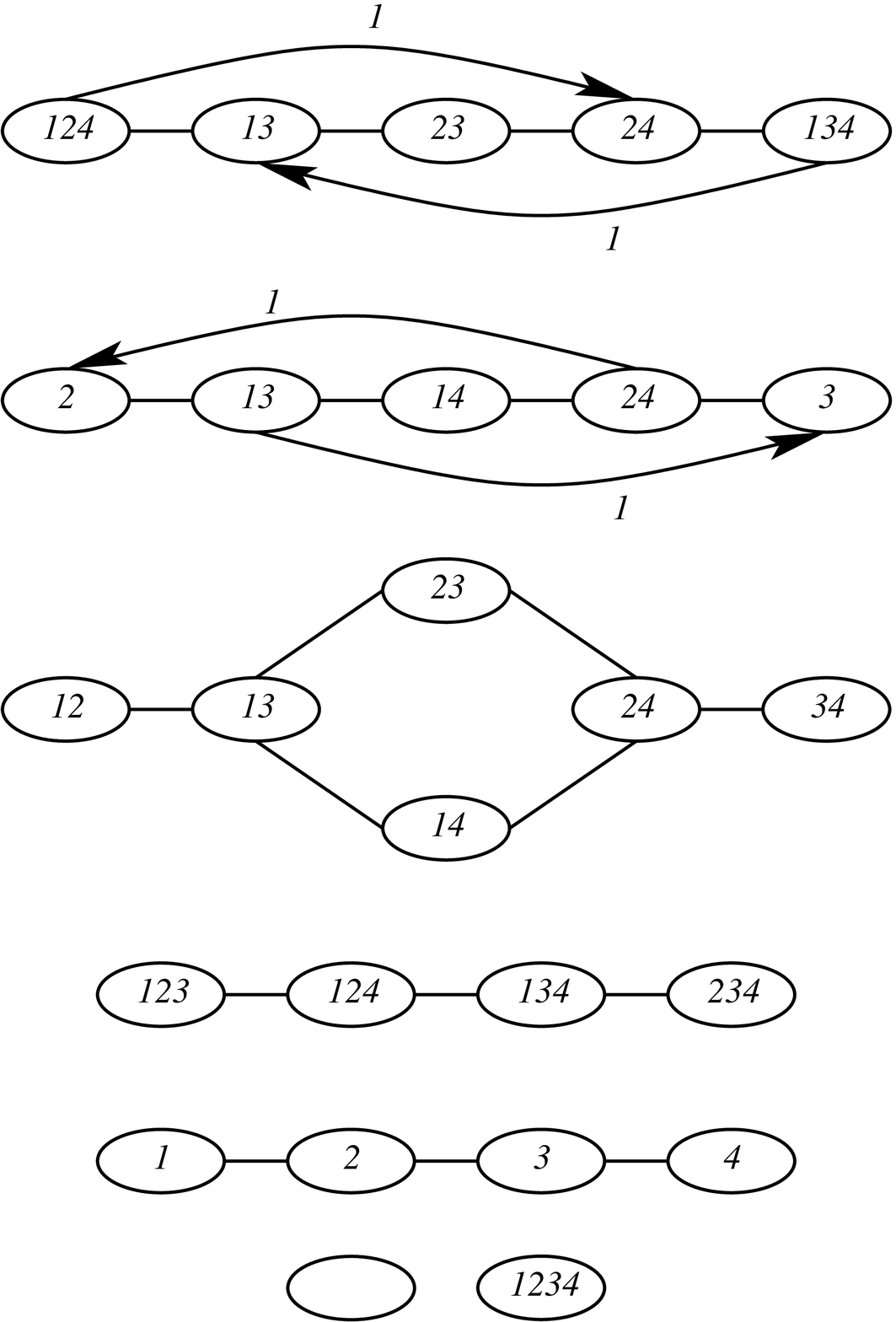}\quad &\quad \includegraphics[height=7.5cm]{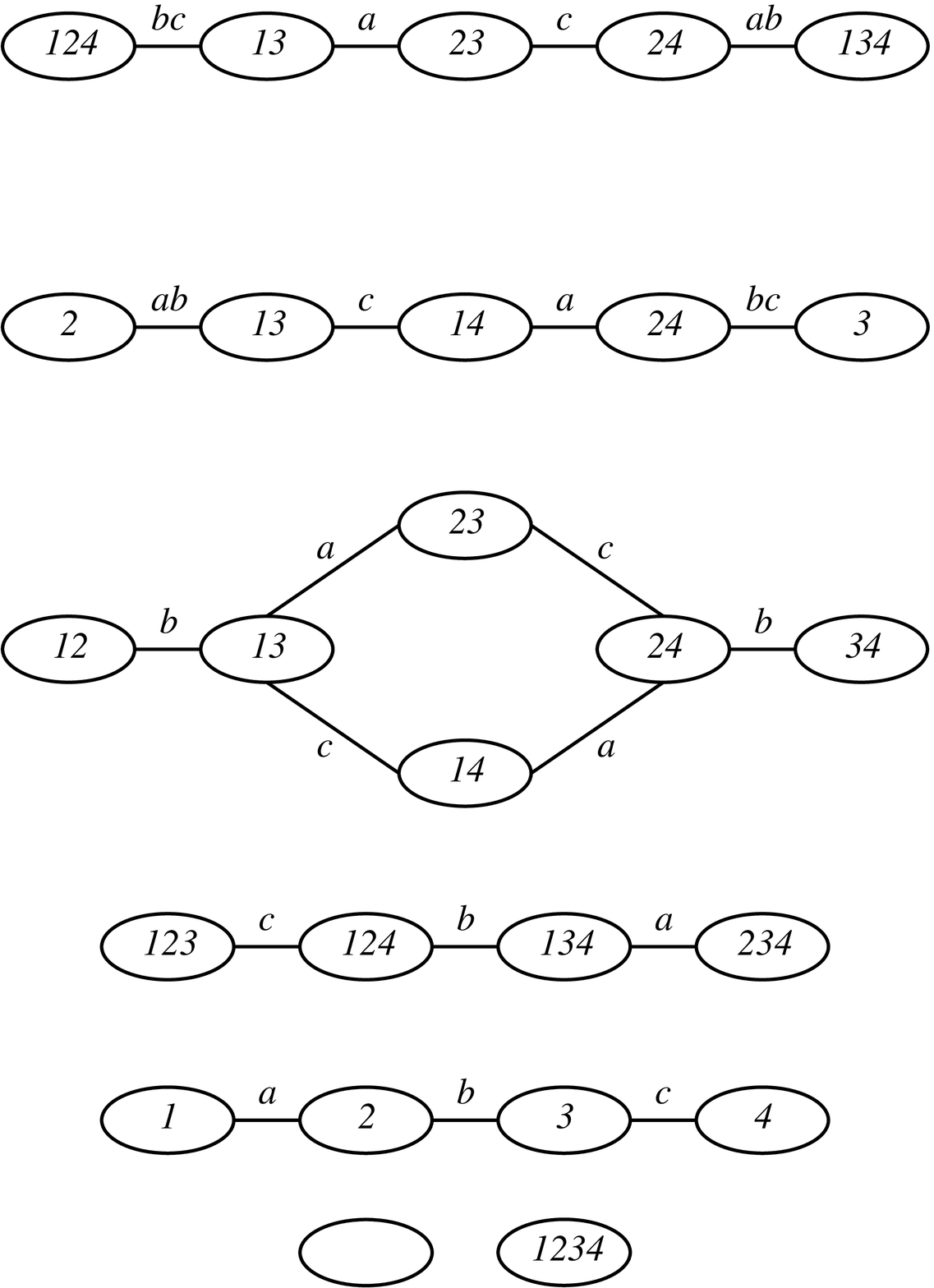}
\end{tabular}
\end{figure}
\end{ex}

The \emph{simple part} of a molecule is the graph formed by erasing all the arcs. We usually view it as an undirected graph. A \emph{morphism of molecules} $\varphi:M\to N$ is a map between the 
vertex sets which
\begin{enumerate}
\item is a graph morphism of the simple parts,
\item preserves $\tau$-invariants.
\end{enumerate}
Notice that a morphism does not need to respect arcs.

\subsection{Restriction}

Let $J\subseteq S$ and let $W_J$ be the corresponding parabolic subgroup.

Let  $M = (V, m, \tau)$ be a $W$-molecular graph. The \emph{$W_J$-restriction} of $M$ is $N=(V, m', \tau')$, with
\begin{enumerate}
\item for all $v\in V$, $\tau'(v) = \tau(v)\cap J$,
\item for all $u,v\in V$, 
$$m'(u, v)=
\begin{cases}
0        , \quad\text{if $\tau'(u)\subseteq\tau'(v)$,}\\
m(u, v), \quad\text{otherwise.}
\end{cases}$$
\end{enumerate} 
The $W_J$-restriction of $M$ is a $W_J$-molecular graph. A {\it $W_J$-submolecule} of $M$ is a $W_J$-molecule (i.e. component 
connected by simple edges) of the $W_J$-restriction of $M$. There is a natural inclusion map of a $W_J$-submolecule into the original molecular 
graph. Sometimes, abusing notation, we refer to the image of this map as a $W_J$-submolecule. The sense in which we use the word should be clear 
from the context.

\begin{rmk}
\label{rmk:lpr restriction}
It is sometimes convenient to think of the local polygon rules in terms of restriction. For LPR2, suppose $k\in\tau(v)\setminus\tau(u)$. Then 
LPR2 holds for an $S$-labeled graph if and only if it holds for the $W_{\{i,j,k\}}$ restriction. Similarly, LPR3 holds for an $S$-labeled graph if 
and only if it holds for the $W_{\{i,j,k,l\}}$ restriction.
\end{rmk}

\subsection{Arc transport}
The following three propositions follow from the local polygon rules. In fact, although we do not need it, for simply laced types they are equivalent to the local polygon rules.

Before proceeding we would like to make a small comment about the bipartition requirement of admissible $S$-labeled graphs. The proposition makes a statement that certain edges have equal weights. Since we are only shown parts of the graph, it is unclear whether the indicated edges break bipartition. However the edges are allowed or disallowed simultaneously. In case the edges are disallowed, the statements of the propositions are trivial. So in the proofs we assume that the edges are, in fact, allowed. 

\begin{prop}
\label{prop:at1}
Suppose $M$ is a $W$-molecular graph. Suppose $(x,y)$ and $(x',y')$ are simple edges that activate the same bond, say $(i,j)$. 
Without loss of generality, $i\in\tau(x)\cap \tau(x')$ and $j\in\tau(y)\cap\tau(y')$. Suppose moreover that there exists $k$ 
such that $k\in\tau(x)\cap\tau(y)$ and $k\notin\tau(x')\cup\tau(y')$. Then $m(x, x') = m(y, y')$. In picture notation (after restriction to the parabolic subgroup generated by $J=\{i,j,k\}$), 
the blue (dashed) edges have the same weight:
$$\includegraphics[width=3cm]{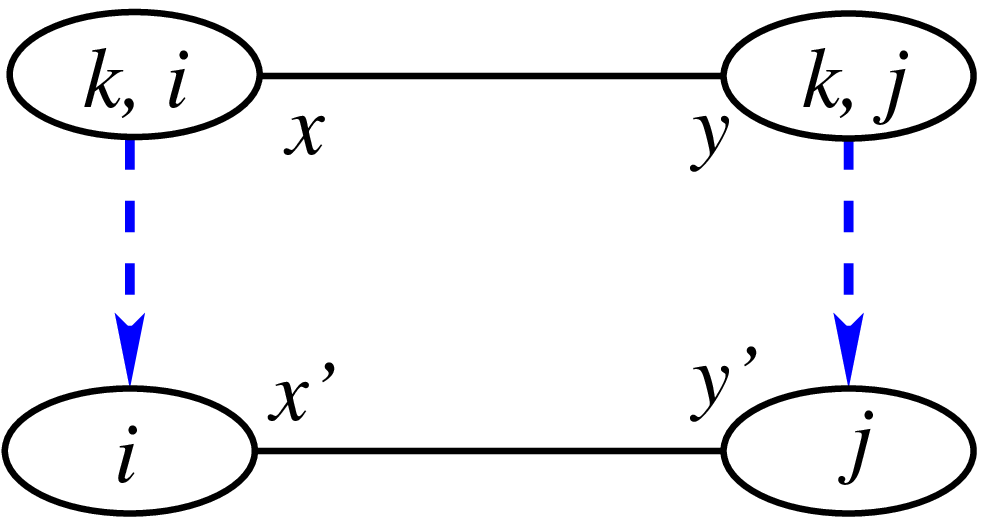}$$
\begin{proof}

There are two evident instances of LPR2, namely
$$N_{ki}^2(M,x, y') = N_{ik}^2(M,x, y'),$$
and
$$N_{kj}^2(M,y, x') = N_{jk}^2(M,y, x').$$

Consider the first of these. Let us analyze when could there be other possible alternating paths besides the ones pictured.
First look at alternating paths of type $(k,i)$. They must pass through a vertex $z\neq y$ with $k\in\tau(z)$ and $i\notin\tau(z)$. 
Since $z\neq y$ and $(x,z)$ is an edge, BR tells us that $j\notin\tau(z)$. Now $(z,y')$ must be an edge, so by CR we have that $(j,k)$ is a bond. 
Hence $N_{ki}^2(M,x, y') = m(y, y')$ unless $(j,k)$ is a bond.

Now look at alternating paths of type $(i,k)$. They must pass through a vertex $z\neq x'$ with $i\in\tau(z)$ and $k\notin\tau(z)$. 
Since $z\neq x'$ and $(z,y')$ is and edge, BR tells us that $j\in\tau(z)$. Now $(x,z)$ must be an edge, so by CR we have that $(j,k)$ is a bond. 
Hence $N_{ik}^2(M,x, y') = m(x, x')$ unless $(j,k)$ is a bond.

Thus the first instance of LPR2 give the desired result unless $(j,k)$ is a bond. By the same argument with $i$ and $j$ switched,
the second instance of LPR2 give the desired result unless $(i,k)$ is a bond. But the Coxeter graph cannot contain triangles, so 
at least one of $(i,k)$ and $(j,k)$ cannot be a bond. This finishes the proof.
\end{proof}
\end{prop}

\begin{prop}
\label{prop:at2}
Suppose $M$ is an $A_3$-molecular graph which contains the simple edges of one of the two pictures:
$$\includegraphics[width=4.5cm]{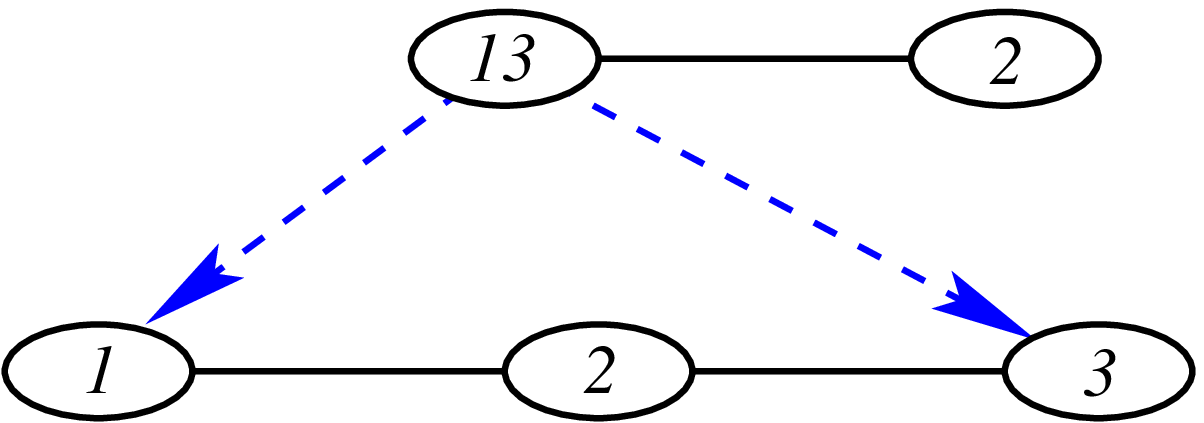},\qquad \includegraphics[width=4.5cm]{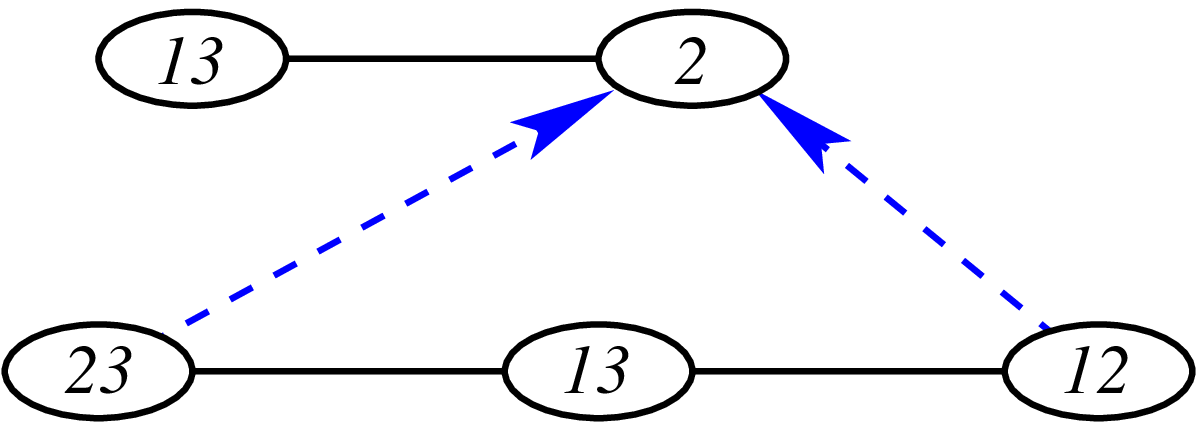}.$$
Then the weights of the two blue (dashed) edges are equal.
\begin{proof}
Apply the only evident instance of LPR2. By Remark \ref{rmk:lprpaths} we are seeing all the possible paths involved, so the desired equality follows. 
\end{proof}
\end{prop}

\begin{prop}
\label{prop:at3}
Suppose $M$ is an $A_4$-molecular graph which contains the simple edges of one of the two pictures:
$$\includegraphics[width=7cm]{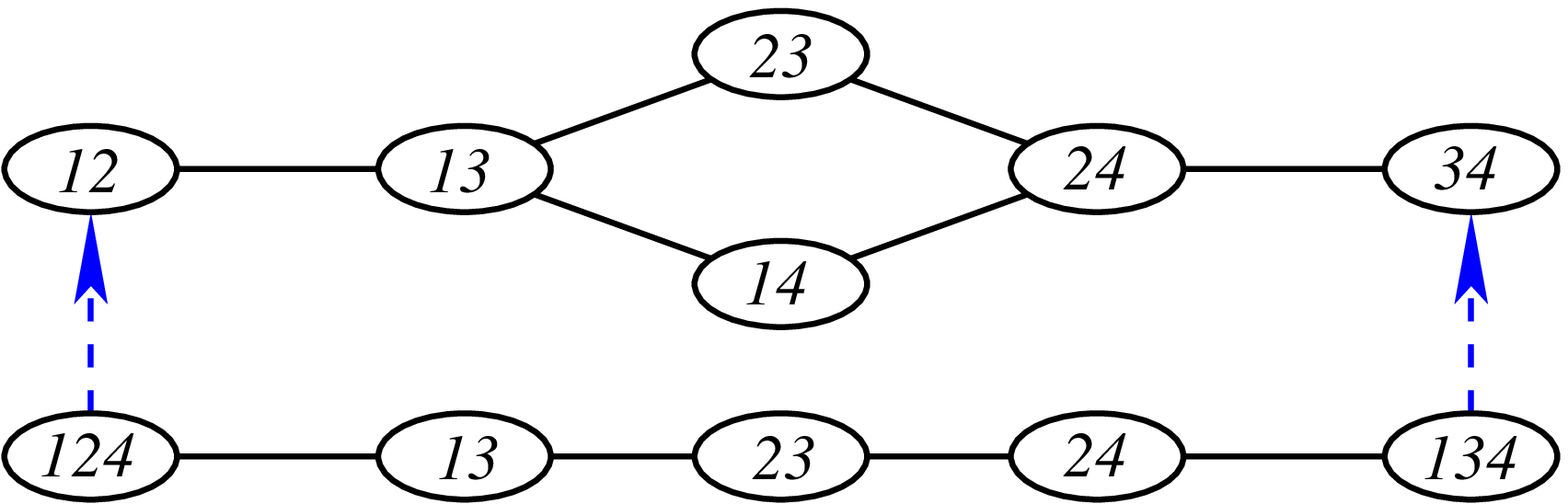},\qquad \includegraphics[width=7cm]{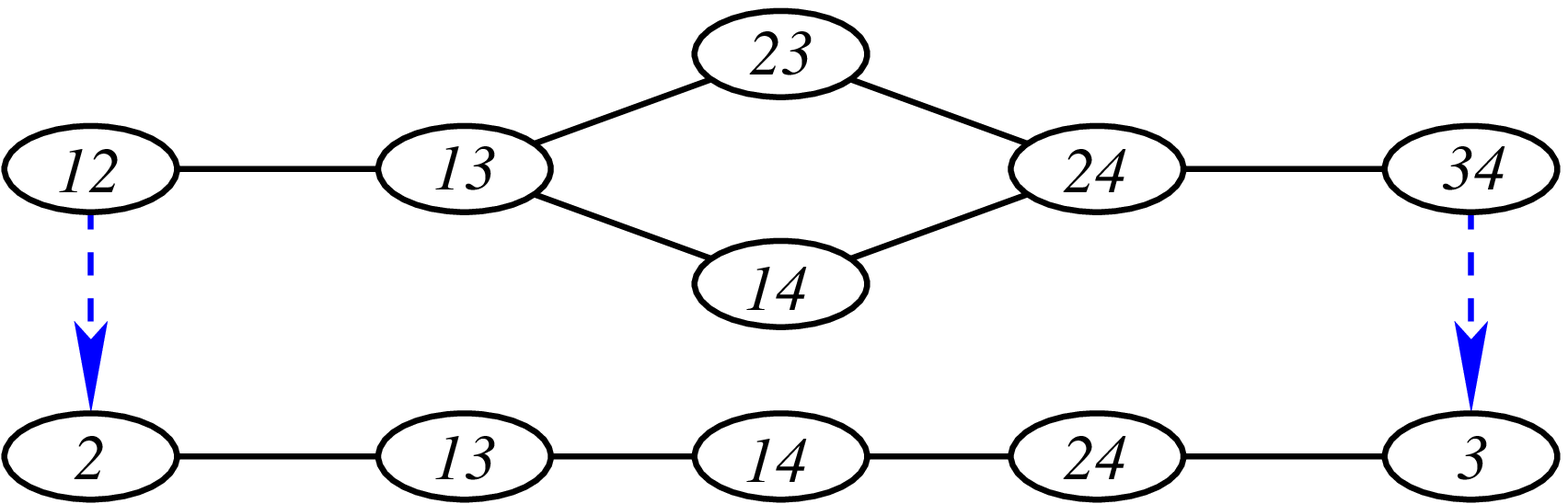}.$$
Then the weights of the two blue (dashed) edges are equal.
\begin{proof}

We will just treat the left picture; the right one is done in the same way. 

Label some of the vertices as follows:
$$\includegraphics[width=7cm]{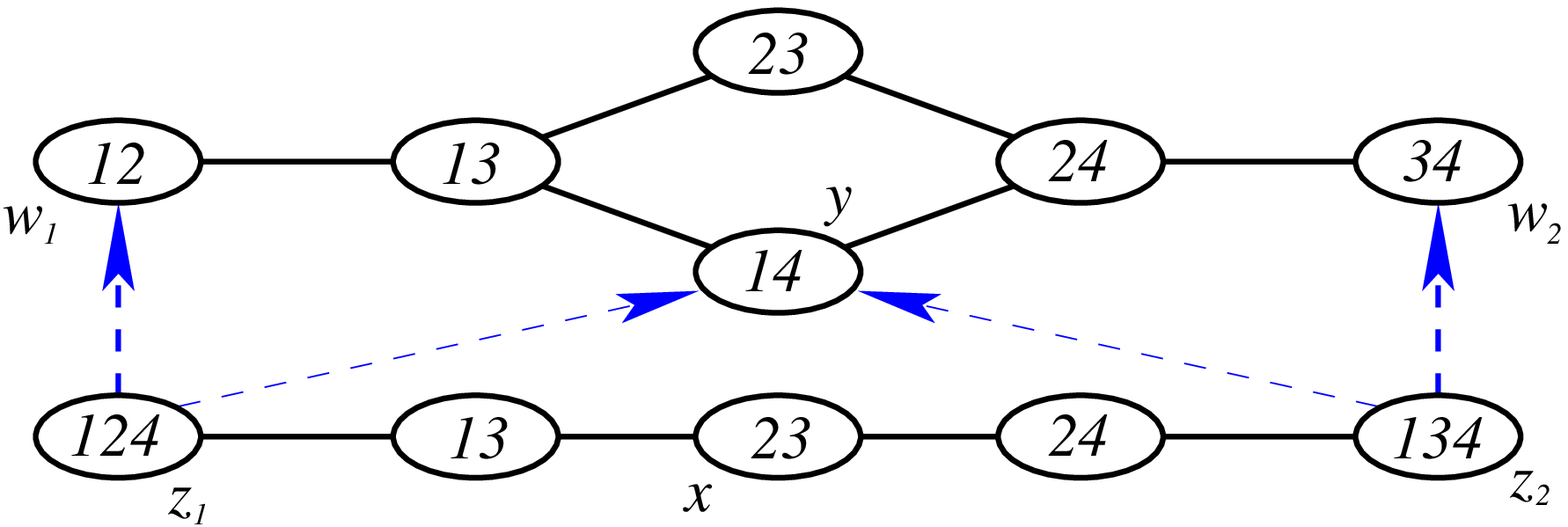}.$$
Applying LPR3 with regard to paths from $x$ to $y$ gives $m(z_1,y) = m(z_2,y)$; as in the proof of Proposition \ref{prop:at2} we can see all the possible paths. 
Now Proposition \ref{prop:at2} itself gives us that $m(z_1, y) = m(z_1, w_1)$ (to apply the proposition we restrict $M$ to the parabolic subgroup generated by $J=\{2,3,4\}$), 
and $m(z_2, y) = m(z_2, w_2)$. This finishes the proof.
\end{proof}
\end{prop}

\section{Dual equivalence graphs}
\label{sec:degs}
This section summarizes the relevant definitions and results of \cite{assaf_degs}. The results are restated to make the similarity with the $W$-molecule world more apparent. 

Fix $n\in\Z^{>0}$. Let $(W,S)$ be a Coxeter system of type $A_{n}$. Identify $S$ in a natural way with $\{1,\dots, n\}$. Define $a_i$ to be the bond $(i,i+1)$. Then $B:=\{a_1,\dots, a_{n-1}\}$ is the set of all bonds. For examples with small $n$ we will use the notation $a,b,c,\dots$ instead of $a_1, a_2, a_3,\dots$.

\begin{defn}
A \emph{signed colored graph of type $n+1$} is a tuple $(V, E, \tau, \beta)$, where $(V,E)$ is a finite undirected simple graph, $\tau:V\to 2^{S}$, and $\beta:E\to 2^{B}$.
\end{defn}

Denote by $E_i$ the set of edges with label $i$ (i.e. such that the corresponding value of\textbf{} $\beta$ contains $i$); we call these \emph{$i$-colored} edges. This is a slight reindexing from Assaf's original definition; in the original $E_i$ was the set of edges whose label contains $i-1$.

We start by constructing a family, indexed by partitions, of signed colored graphs. 

\subsection{``Standard'' dual equivalence graphs}
\label{sec:std degs}

Let $\lambda$ be a partition of $n+1$. 
Let $SYT(\lambda)$ be the set of standard Young tableaux of shape $\lambda$. Using the English convention for tableaux, the left-descent set of a tableau $T$ is 
$$\tau(T) := \{1\leqslant i\leqslant n: \text{$i$ is located in a higher row than $i+1$ in $T$}\}.$$
The set of vertices of our graph is $V:=SYT(\lambda)$ (see Example \ref{ex:degs}).

By a diagonal of a tableau we mean a $NW-SE$ diagonal.
A \emph{dual Knuth move} is the exchange of $i$ and $i+1$ in a standard tableau, provided that either $i-1$ or $i+2$ lies (necessarily strictly) 
between the diagonals containing $i$ and $i+1$. This corresponds to dual Knuth moves on the symmetric group via, for example, the ``content reading word'' 
(reading each diagonal from top to bottom, and concatenating in order of increasing height of the diagonals). 
The dual Knuth moves define the edges of our graph:
$$E:=\{(T,U): \text{$T$ and $U$ are related by a dual Knuth move}\}.$$

A dual Knuth move between tableaux $T$ and $U$ \emph{activates the bond $a_i$} if $i$ lies in precisely one of $\tau(T)$ and $\tau(U)$, and $i+1$ lies 
precisely in the other. Denote this condition by $T\stackrel{a_i}{-} U$.
For $(T, U)\in E$, let
$$\beta(T,U) := \{a_i\in B: T\stackrel{a_i}{-} U\}.$$

The graph $G_\lambda:=(V, E, \tau, \beta)$ is a signed colored graph of type $n+1$.

One can give a slightly more explicit description of activations on tableaux. Notice that $i, i+1$, and $i+2$ have to 
lie on three distinct diagonals in any tableau. We have $T\stackrel{a_i}{-} U$ precisely when $T$ and $U$ differ 
by switching the two of the above entries on the outside diagonals, provided that the middle diagonal does not contain $i+1$.

\begin{ex}
\label{ex:degs}
Here are two standard dual equivalence graphs, corresponding to the shapes $311$ and $32$.
$$\includegraphics[width=16cm]{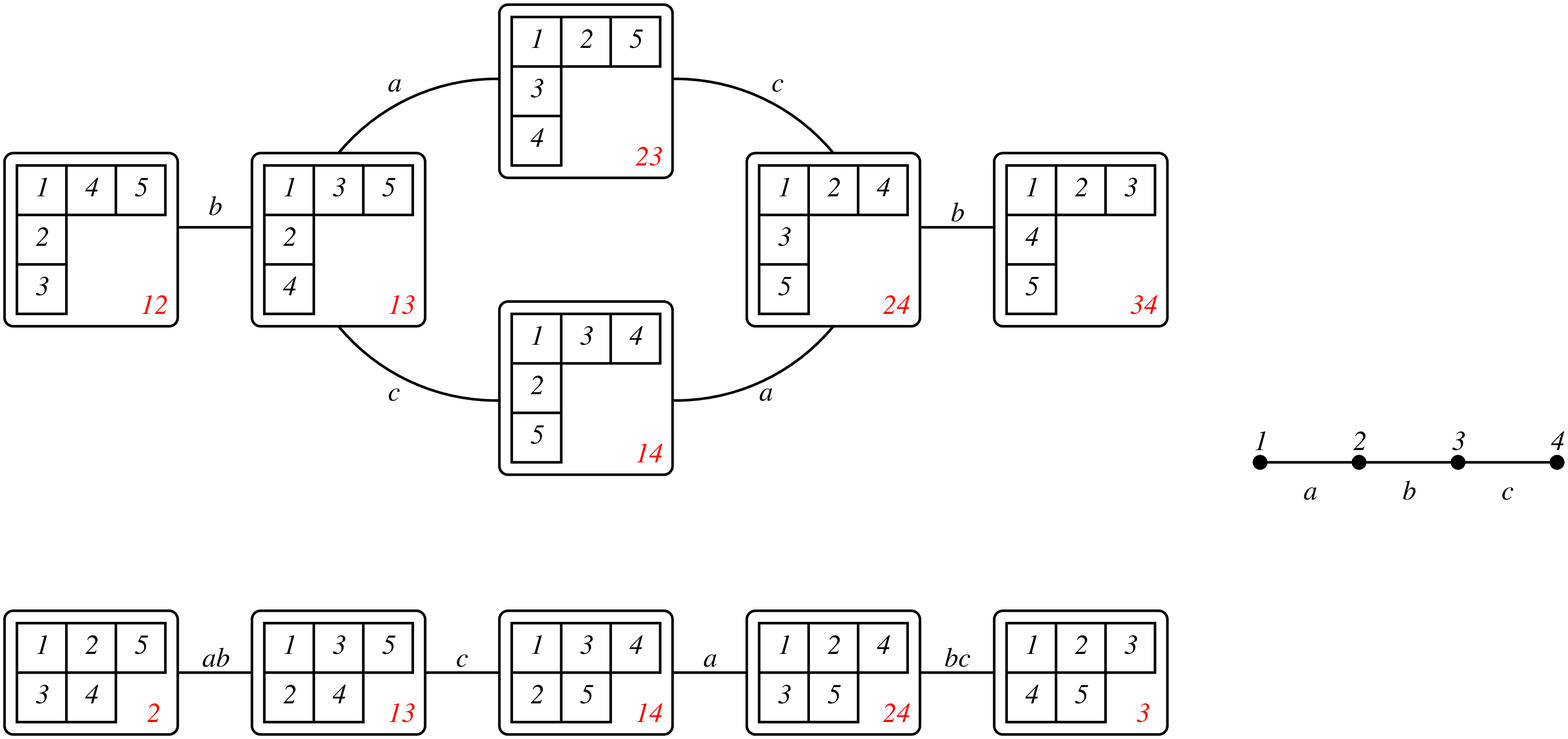}$$
The values of $\tau$ are shown in red in the lower right-hand corner of each vertex.
\end{ex}

\subsection{Axiomatics}
Now we review Assaf's axiomatization of the above construction. 

A vertex $w$ of a signed colored graph is said to \emph{admit an $i$-neighbor} if precisely one of $i$ and $i+1$ lies in $\tau(w)$.

\begin{defn}
A \emph{dual equivalence graph of type $n+1$} is a signed colored graph $(V, E, \tau, \beta)$ such that for any $1\leqslant i < n	$:
\begin{enumerate}
\item For $w\in V$, $w$ admits an $i$-neighbor if and only if there exists $x\in V$ which is connected to $w$ by an edge of color $i$. In this case $x$ must be unique.
\item Suppose $(w,x)$ is an $i$-colored edge. Then $i\in\tau(w)$ iff $i\notin\tau(x)$, $i+1\in\tau(w)$ iff $i+1\notin\tau(x)$, and if $h<i-1$ or $h>i+2$ then 
$h\in\tau(w)$ iff $h\in\tau(x)$.

In other words, going along an $i$-colored edge switches $i$ and $i+1$ in the $\tau$-invariant, and does not affect any labels except $i-1, i, i+1$, and $i+2$.   
\item Suppose $(w,x)$ is an $i$-colored edge. If $i-1\in\tau(w)\Delta\tau(x)$ then ($i-1\in\tau(w)$ iff $i+1\in\tau(w)$), where $\Delta$ is the symmetric difference. If $i+2\in\tau(w)\Delta\tau(x)$ then 
($i+2\in\tau(w)$ iff $i\in\tau(w)$).
\item If $i<n-2$, consider the subgraph on all the vertices and $i$- and $(i+1)$-colored edges. Each of its connected components has the form:
$$\includegraphics[height=.4cm]{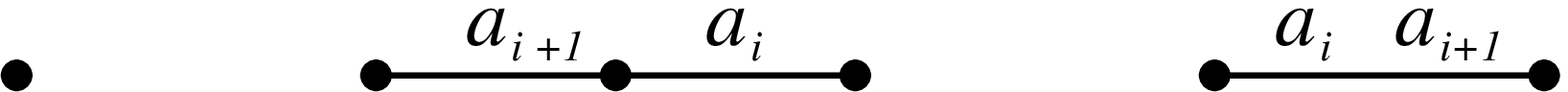}.$$
If $i<n-3$, consider the subgraph on all the vertices and $i$-, $(i+1)$-, and $(i+2)$-colored edges. Each of its connected components has the form:
$$\includegraphics[height=2.6cm]{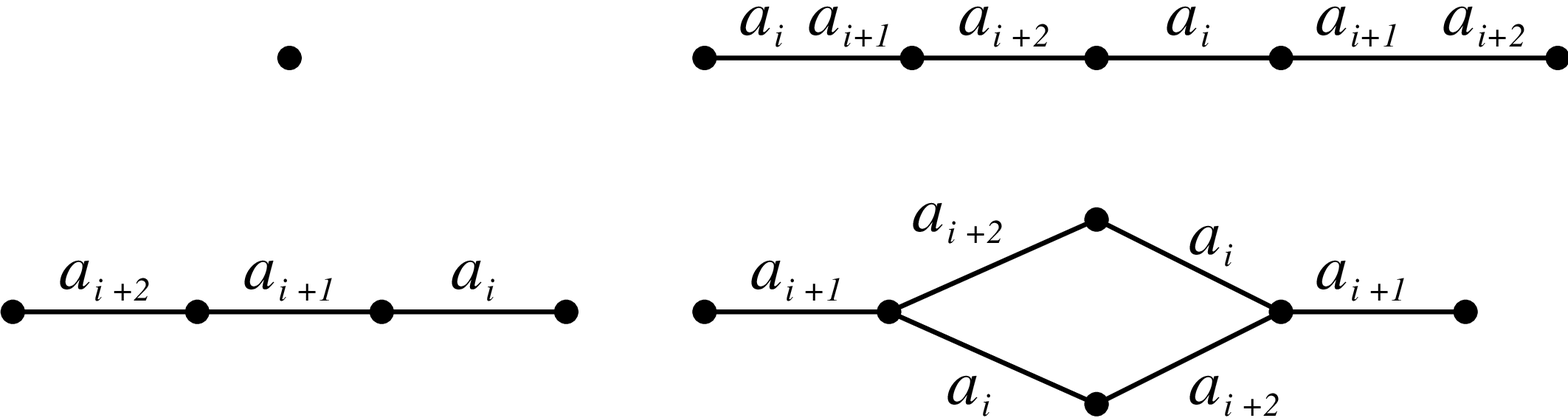}.$$
\item Suppose $(w,x)\in E_i, (x,y)\in E_j$, and $\abs{i-j}\geqslant 3$. Then there exists $v\in V$ such that  $(w,v)\in E_j, (v,y)\in E_i$.
\item Consider a connected component of the subgraph on all the vertices and edges of colors $\leqslant i$. If we erase all the $i$-colored edges it breaks down into several components. Any two of these were connected by an $i$-colored edge.
\end{enumerate}
\end{defn}

Examples of $A_4$ dual equivalence graphs can be found on the right of Figure \ref{fig:a4molsdegs}.

A \emph{morphism} of signed colored graphs is a map on vertex sets which preserves $\tau$ and $\beta$. 

\begin{prop}
The graph $G_\lambda$ is a dual equivalence graph. Moreover, $\{G_\lambda\}_\lambda$ is a complete collection of isomorphism class representatives of dual equivalence graphs. 
\begin{proof}
The references are to \cite{assaf_degs}. The first statement is Proposition 3.5. The second is a combination of Theorem 3.9 and Proposition 3.11.
\end{proof}
\end{prop}

\begin{rmk}
There is some redundancy in the definition as presented. Namely, $\beta$ can be calculated from $\tau$: an edge $(u,v)$ is $i$-colored if and only 
if $i$ lies in precisely one of $\tau(u)$ and $\tau(v)$, while $i+1$ lies in precisely the other. Assaf needed a slightly more general definition, 
so $\beta$ was not redundant. We think of $\beta$ as a piece of data about a dual equivalence graph, and keep it as 
part of the definition to be consistent with the original. 
\end{rmk}

A \emph{weak dual equivalence graph} is a signed colored graph satisfying $1-5$ of the above.

\subsection{Restriction}

Suppose $G$ is a signed colored graph of type $n+1$. For $0\leqslant k<n+1$, a \emph{$(k+1)$-restriction} of $G$ consists of the same vertex 
set $V$, the $\tau$ function post-composed with intersection with $\{1, \dots, k\}$, and the $\beta$ function post-composed with restriction to 
$\{a_1, \dots, a_{k-1}\}$. The $(k+1)$-restriction of $G$ is a signed colored graph of type $k+1$. The property of being a (weak) dual equivalence 
graph is preserved by restriction. By a \emph{$(k+1)$-component of $G$} we mean either the connected component of the restriction, or the induced 
subgraph of $G$ on vertices corresponding to such connected component. It should be clear from the context which of these we are talking about.

The $n$-components of $G_\lambda$ are obtained by fixing the position of $n+1$ in the tableau. Such a component is isomorphic to $G_\mu$, 
where $\mu$ is formed from $\lambda$ by erasing the outer corner which contained $n+1$. On the above examples this looks as follows:

\begin{ex}
\label{ex:degs_restr}
$$\includegraphics[width=\textwidth]{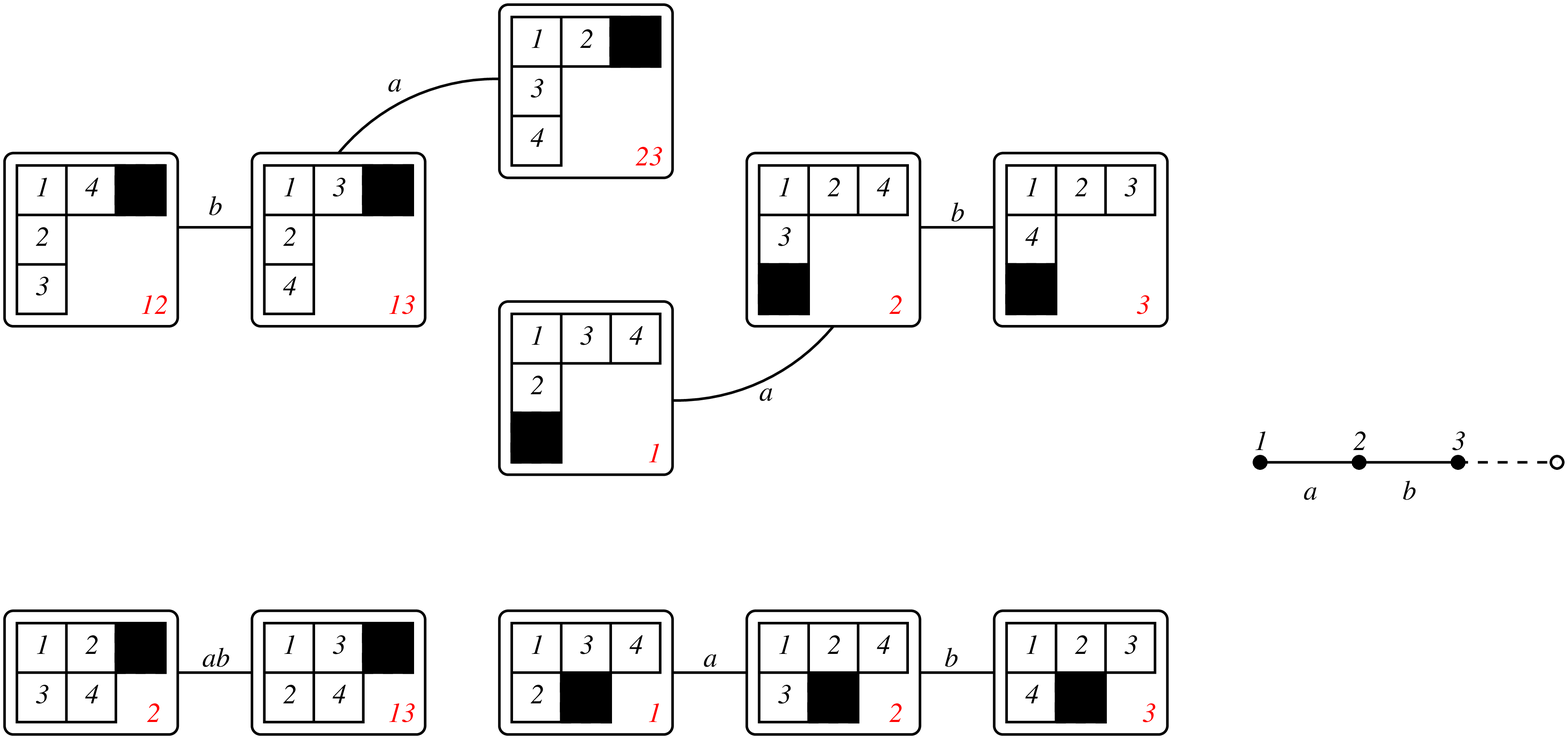}$$
\end{ex}

The condition of being a weak dual equivalence graph is already quite powerful. The following lemma is relevant to us. It essentially says that 
a weak dual equivalence graph with a nice restriction property is necessarily a cover of a dual equivalence graph.

\begin{lemma}
\label{sami:surjective map}
Suppose $G$ is a weak dual equivalence graph of type $n+1$. Suppose moreover that each $n$-component is a dual equivalence graph. 
Then there is a surjective morphism $\varphi: G\to G_\lambda$ for some partition $\lambda$ of $n+1$, 
which restricts to an isomorphism on the $n$-components.

Let $C\cong G_\mu$ be an $n$-component. Then for any partition $\nu\neq\mu$ of $n$ with $\nu\subset\lambda$, there exists a unique $n$-component 
$D$ with $\varphi(D)=G_\nu$ which is connected to $C$ by an $E_{n-1}$ edge. Also, two $n$-components which are isomorphic to $G_\mu$ are 
not connected by an $E_{n-1}$ edge.

\begin{proof}
The references are again to \cite{assaf_degs}. The existence of the morphism is shown in Theorem 3.14. Its surjectivity follows by Remark 3.8. The 
fact that it restricts to an isomorphism on then $n$-components follows from the proof of Theorem 3.14. The covering properties from the second 
paragraph are shown in Corollary 3.15, though the last one is not explicitly mentioned. 
\end{proof}
\end{lemma}

\subsection{Molecules and dual equivalence graphs}

\begin{prop}
The simple part of an $A_n$-molecule, with the corresponding $\tau$ function and edges labeled by activated bonds, is a weak dual equivalence graph.

\begin{proof}
Axioms (1), (2), (3) follow directly from SR, BR, and CR. Axiom (4) was demonstrated in Examples \ref{ex:a3mols} and \ref{ex:a4mols}. Axiom (5) is a weaker version of the local polygon rule.
\end{proof}
\end{prop}

Consider the graph $G_\lambda$. It is clear that (viewed as a directed graph with edge weights of 1) it is an admissible $S$-labeled graph for the $A_n$ root system. It is well known that it forms the simple part of an $A_n$-molecule  (the left Kazhdan-Lusztig cell) which we call $\overline{G_\lambda}$. 

\begin{defn} 
An $A_n$-molecule is \emph{Kazhdan-Lusztig} if it is isomorphic to $\overline{G_\lambda}$, i.e. if its simple part is a dual equivalence graph.
\end{defn}

\begin{rmk}
\label{rmk:submolecules of G_lambda}
We can explicitly describe the simple edges of the parabolic restriction of $\overline{G_\lambda}$. Let $J=\{j_1, \dots, j_k\}$. Then the relevant tableau entries are $J':= \{j_1,j_1+1,j_2,j_2+1, \dots, j_k+1\}$. The simple edges of the $W_J$-restriction of $G_\lambda$ are dual Knuth moves that exchange two entries of $J'$ provided the ``witness'' between them is also in $J'$.
\end{rmk}

\section{Main theorem}
\label{sec:thm}

In this section we show that any $A_n$-molecule is Kazhdan-Lusztig. The proof will proceed by induction on $n$, so the preliminary results 
will start with an $A_n$-molecule whose $A_{n-1}$-submolecules are Kazhdan-Lusztig.

The first of these results states that if two such $A_{n-1}$-submolecules are connected by a simple edge, then the connected $A_{n-2}$-submolecules are isomorphic and there is a ``cabling'' of edges (possibly arcs) of weight 1 between them:
$$\includegraphics[width=5cm]{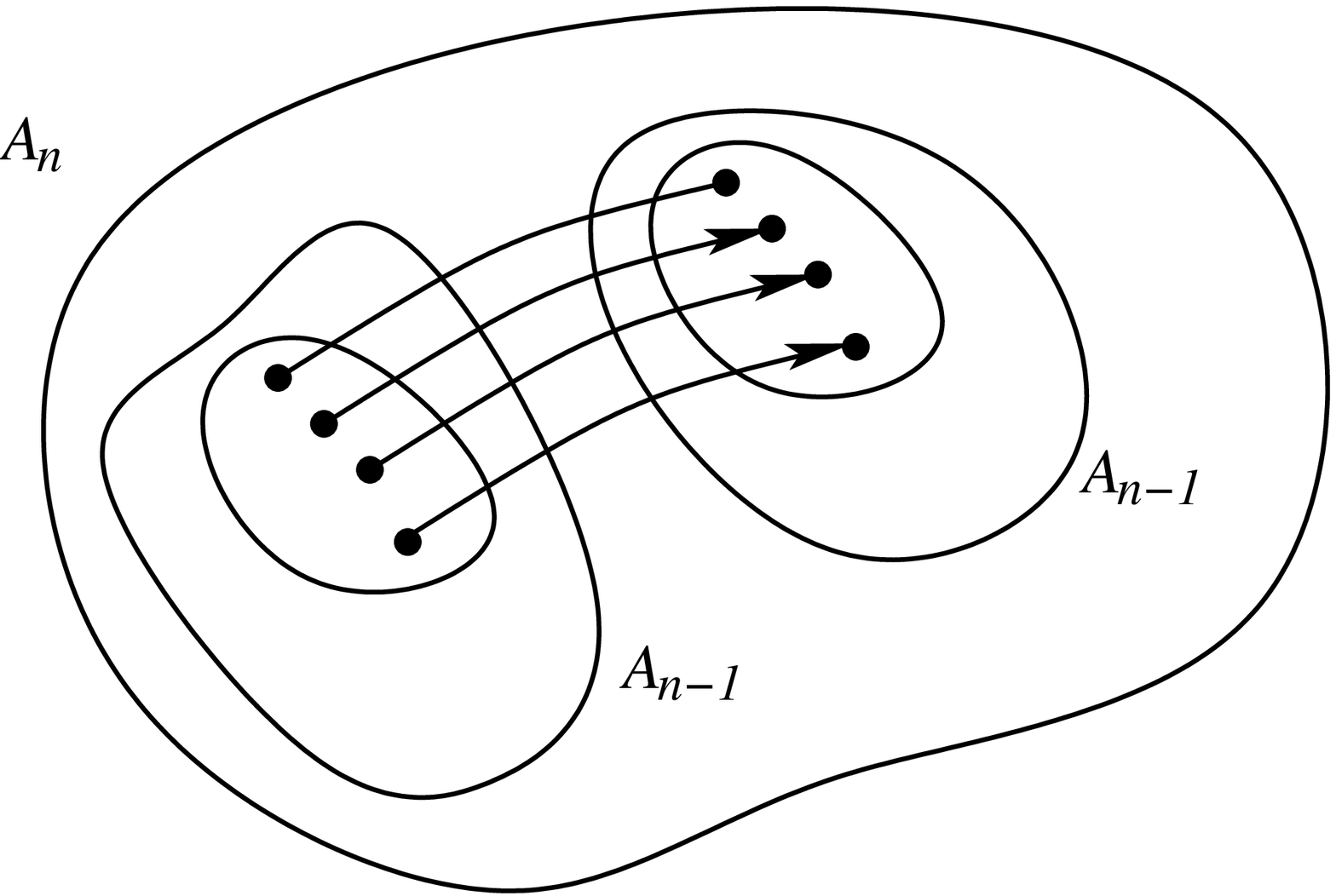}.$$

\begin{lemma}
\label{lemma:cabling}
Let $M$ be an $A_n$-molecule whose $A_{n-1}$-submolecules are Kazhdan-Lusztig. Suppose $A$ and $B$ are two such 
submolecules which are joined by a simple edge (in $M$), namely there exist $x\in A$, $y\in B$ such that the edge 
$(x, y)$ is simple. Let $A'$ (resp. $B'$) be the $A_{n-2}$-submolecule of $M$ containing $x$ (resp. $y$). Then there is 
an isomorphism $\psi$ between $A'$ and $B'$ such that $\psi(x)=y$. Moreover, if $n\in\tau(x)$ then 
$m(z,\psi(z)) = 1$ for all $z\in A'$.
\begin{proof}
By Lemma \ref{sami:surjective map} we know that there is a surjective morphism $\varphi:M\to \overline{G_\lambda}$ for some $\lambda$. Then 
$\varphi(A)\cong \overline{G_\mu}$ and $\varphi(B)\cong \overline{G_\nu}$, for some $\mu, \nu$ which are formed from $\lambda$ by erasing an outer corner 
(these outer corners must be different since no two molecules corresponding to the same shape may be connected; Lemma 
\ref{sami:surjective map}). 

Let $T = \varphi(x), U=\varphi(y)$. Thus $T$ has $n+1$ in position $\lambda\setminus\mu$ and $U$ has 
$n+1$ in position $\lambda\setminus\nu$. Now there is a simple edge between $T$ and $U$, i.e. one is obtained from the other by 
a Knuth move. The only Knuth move in $A_n$ which moves $n+1$ is one that exchanges $n$ and $n+1$, in the presence of $n-1$ between 
them. Hence $T$ has $n$ in position $\lambda\setminus\nu$ and $U$ has $n$ in position $\lambda\setminus\mu$. Hence the $A_{n-2}$-molecule 
containing $T$ has standard tableaux on $\lambda\setminus (\mu\cup \nu)$ as vertices, and Knuth moves between them as edges. The same is 
true for the $A_{n-2}$ molecule containing $U$. Thus the two molecules are isomorphic and the isomorphism is given by switching $n$ and $n+1$. 
Now we can use $\varphi$ to lift it up to an isomorphism $\psi:A'\to B'$.

Suppose $n\in\tau(x)$. Then $n\notin\tau(y)$. So for any $x'\in A'$,  we have $n\in\tau(x')$, and similarly for any $y'\in B'$, $n\notin\tau(y')$. 
Then repeated application of Proposition \ref{prop:at1} shows that the weight of the edge between $z\in A'$ and $\psi(z)$ is the same as between 
$x$ and $y$, namely $1$.
\end{proof}
\end{lemma}

The second preliminary result shows that if, out of three $A_{n-1}$-submolecules, two pairs (satisfying some conditions) are connected by 
simple edges, then the third pair is also connected by a simple edge: 
$$\includegraphics[width=5cm]{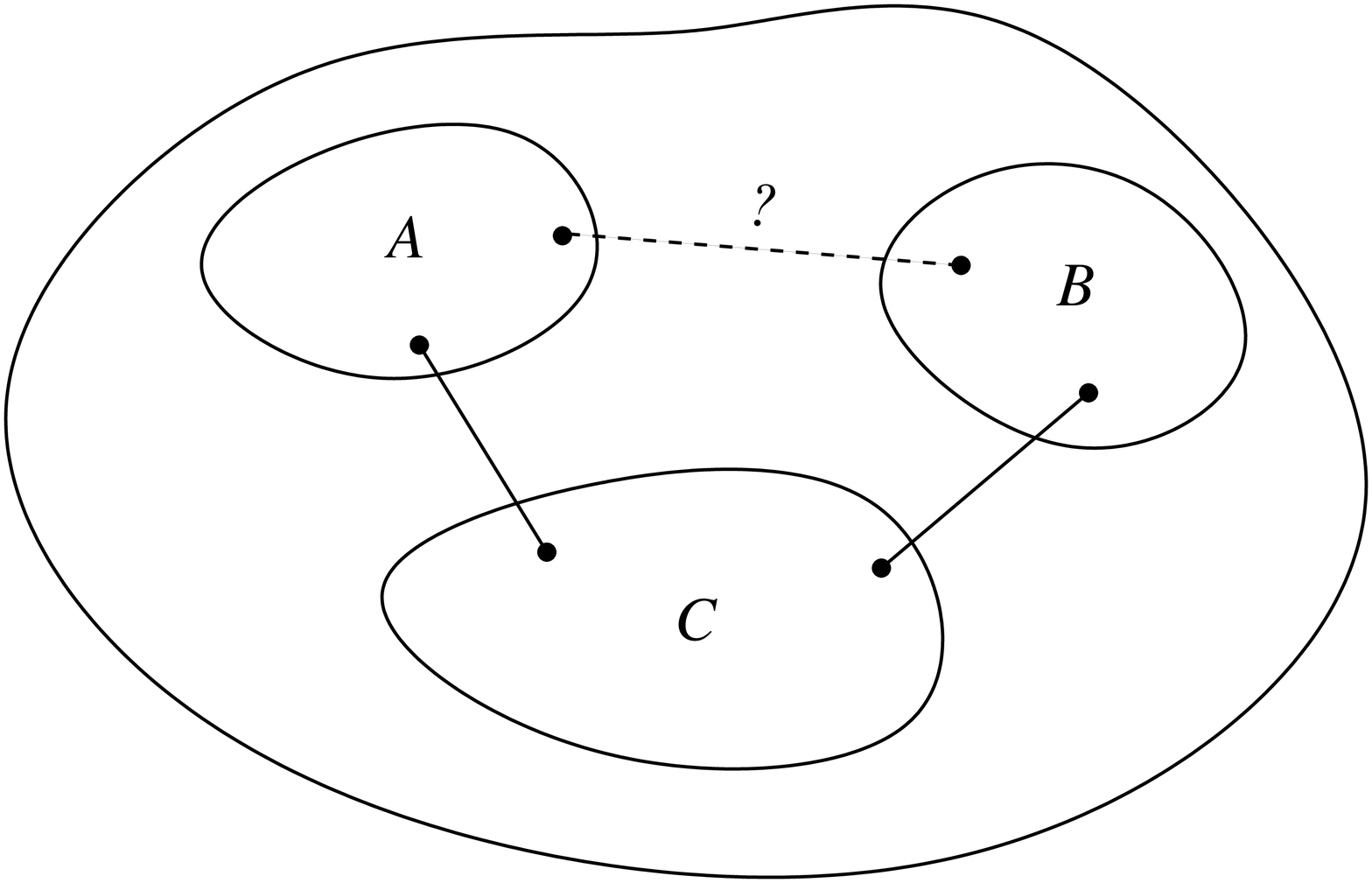}.$$
The conditions will later be removed to show that any two $A_{n-1}$-submolecules of an $A_n$-molecule are connected by a simple edge.

\begin{lemma}
\label{lemma:shortening}
Let $M$ be an $A_n$-molecule whose $A_{n-1}$-submolecules are Kazhdan-Lusztig. By Lemma \ref{sami:surjective map}, there is a surjective 
morphism $\varphi:M\to \overline{G_\lambda}$ for some partition $\lambda$ of $n+1$. Let $A,B,C$ be $A_{n-1}$-submolecules of $M$ such that 
$A$ and $B$ are both connected to $C$ by simple edges. 
Then $A\cong \overline{G_\mu}$, $B\cong \overline{G_\nu}$, $C\cong \overline{G_\eta}$, for some partitions formed by deleting
outer corners of $\lambda$. The three partitions have to be different by Lemma \ref{sami:surjective map}. Suppose moreover that the deleted 
corner for $\eta$ was the highest of the three, namely:

$$\bigskip\scalebox{.3}{\begin{picture}(0,0)%
\includegraphics{submolecule_partitions.pstex}%
\end{picture}%
\setlength{\unitlength}{3947sp}%
\begingroup\makeatletter\ifx\SetFigFont\undefined%
\gdef\SetFigFont#1#2#3#4#5{%
  \reset@font\fontsize{#1}{#2pt}%
  \fontfamily{#3}\fontseries{#4}\fontshape{#5}%
  \selectfont}%
\fi\endgroup%
\begin{picture}(13266,2766)(118,-2044)
\put(10276,-61){\makebox(0,0)[lb]{\smash{{\SetFigFont{29}{34.8}{\rmdefault}{\mddefault}{\itdefault}{\color[rgb]{0,0,0}$\eta$}%
}}}}
\put(676,-61){\makebox(0,0)[lb]{\smash{{\SetFigFont{29}{34.8}{\rmdefault}{\mddefault}{\itdefault}{\color[rgb]{0,0,0}$\mu$}%
}}}}
\put(5476,-61){\makebox(0,0)[lb]{\smash{{\SetFigFont{29}{34.8}{\rmdefault}{\mddefault}{\itdefault}{\color[rgb]{0,0,0}$\nu$}%
}}}}
\end{picture}%
}$$

Then $A$ and $B$ are connected by a simple edge. 
\begin{proof}
Notice that the role of $A$ and $B$ is symmetric, so without loss of generality we may assume that the deleted corner for $\mu$ was the 
lowest of the three.

To prove the existence of an edge between $A$ and $B$, we will choose a simple edge of $C$ and show, using arc transport rules, that its weight 
is equal to the weight of an edge between a vertex of $A$ and a vertex of $B$ whose $\tau$ invariants are incomparable. This will show that 
the edge in question is simple.

Consider a simple edge in $\overline{G_\lambda}$ (which happens to lie in the submolecule isomorphic to $\overline{G_\eta}$) 	of the form:
$$\includegraphics[width=7cm]{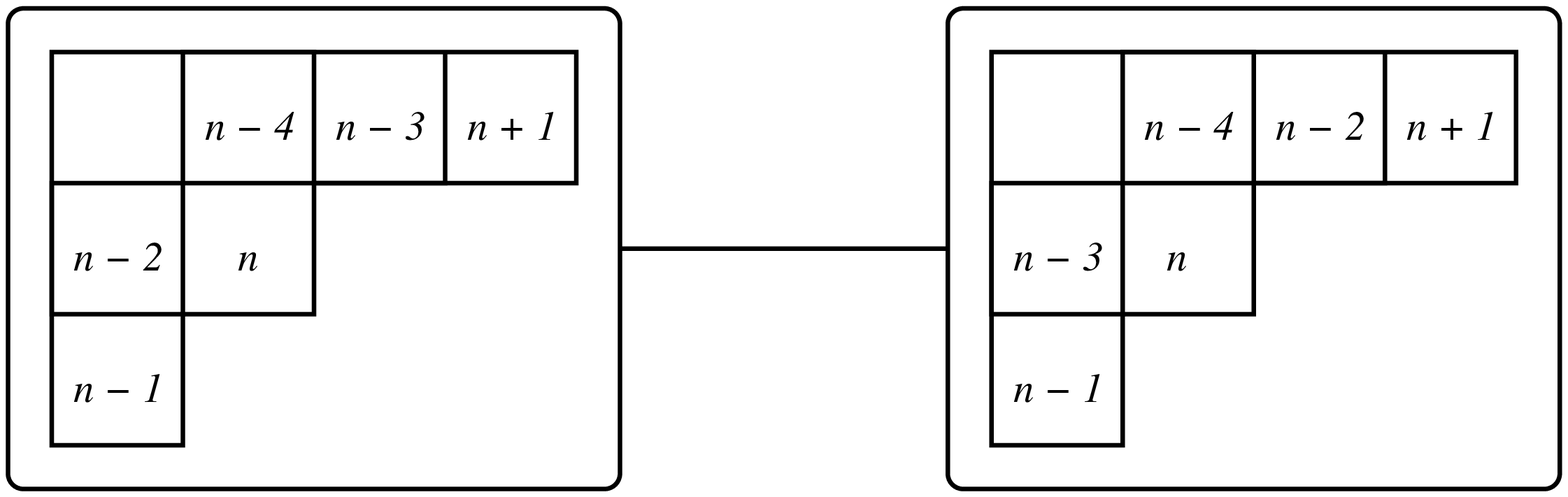}.$$
Let us describe precisely the kind of tableau we are looking for on the left. We want $n+1$ to occupy the cell $\lambda\setminus\eta$, 
$n$ to occupy the cell $\lambda\setminus\nu$, and $n-1$ to occupy the cell $\lambda\setminus\mu$. There exists an outer corner of 
$\mu\cap\nu\cap\eta$ which now lies on a diagonal between $n$ and $n+1$; place $n-3$ there. Place $n-2$ in the outer corner of 
$\mu\cap\nu\cap\eta$ between $n-1$ and $n$. Similarly, place $n-4$ between $n-2$ and $n-3$. Fill in the rest of the tableau in an arbitrary way. The 
two resulting tableaux differ by a Knuth move: one may flip $n-2$ and $n-3$ since $n-4$ is between them. So this is indeed a simple edge in 
$\overline{G_\lambda}$.

Now look at the $A_4$-molecules involved after restricting to the \emph{rightmost} copy of $A_4$ in $A_n$. The restriction corresponds to allowing 
Knuth moves that exchange entries $\geqslant n-3$ provided the ``witness'' between them is also $\geqslant n-3$ (in particular, the original simple 
edge will become directed in the restriction). These are shown in Figure \ref{fig:transport}. The weight of the left blue (directed) edge is $1$ since it 
was a simple edge before the restriction. It is equal to the weight of the right blue (dashed) edge by Proposition \ref{prop:at3}.

\begin{figure}
\caption{Transport along $A_4$ molecules in $\overline{G_\lambda}$.}
\label{fig:transport}
\bigskip\scalebox{.31}{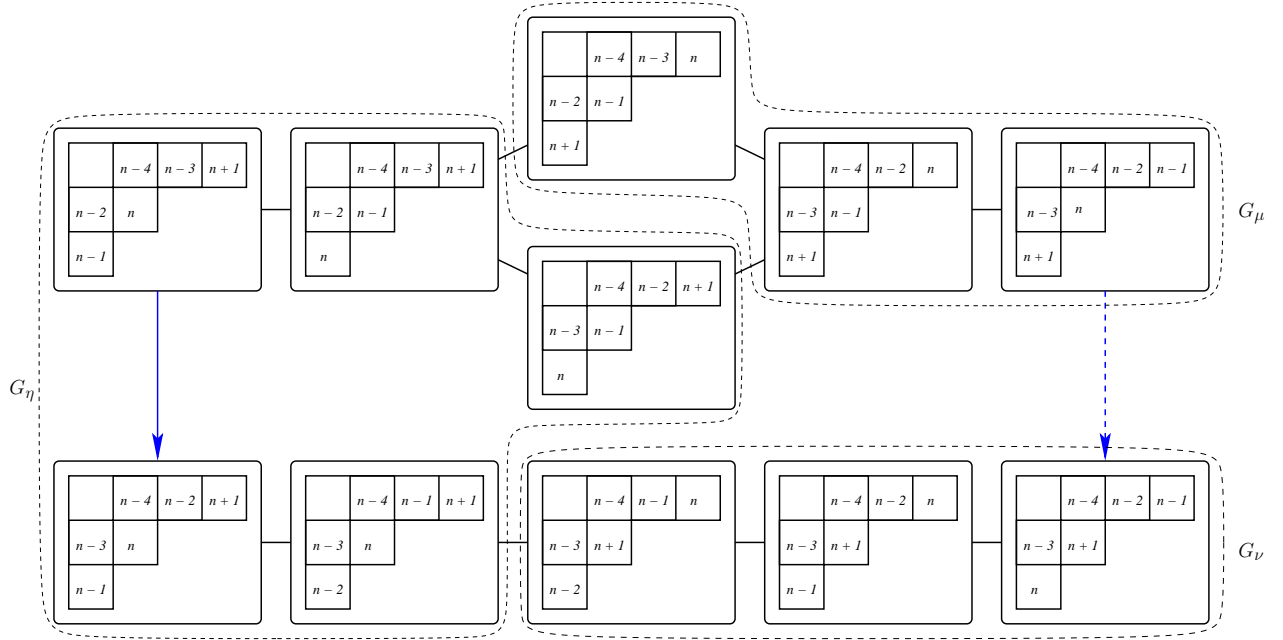}
\end{figure} 

In the original $\overline{G_\lambda}$, before restriction, we may then use the cabling of Lemma \ref{lemma:cabling}, to further transport 
this edge weight as in Figure \ref{fig:transport_ext}.
\begin{figure}
\caption{Transport along a cabling in $\overline{G_\lambda}$.}
\label{fig:transport_ext}
$$\includegraphics[width=6cm]{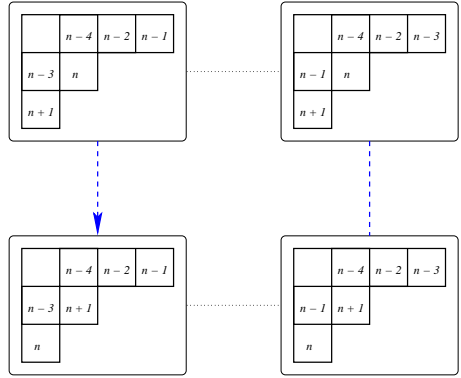}$$
\end{figure} 
Thus we have shown that the weight of the right blue edge in this figure is $1$. 

In $\overline{G_\lambda}$ this is not very interesting since the two tableaux are seen to be related by a Knuth move; let us lift our sequence of 
moves up to $M$. Our original simple edge was located in the submolecule isomorphic to $\overline{G_\eta}$. The preimage under $\varphi$ of 
that simple edge lies in $C$. Now consider the preimages of the two $A_4$-molecules. The preimage of the right end of the molecule on the top will lie in $A$ 
since it is the only molecule isomorphic to $\overline{G_\mu}$ which is connected to $C$ by a simple edge (Lemma \ref{sami:surjective map}). Similarly, the preimage of the right end of the molecule 
on the bottom is in $B$. So there is an edge of weight 1 from $A$ to $B$. The transport along a cabling does not change the $A_{n-1}$-molecules 
involved, however the $\tau$-invariants of the right blue edge in Figure \ref{fig:transport_ext} are manifestly incomparable (one has $n-1$ while 
the other has $n$). This produces a simple edge between $A$ and $B$.
\end{proof}
\end{lemma}

We can now finish the proof of the theorem.

\begin{thm}
Any $A_n$-molecule is Kazhdan-Lusztig.

\begin{proof}
We know that the simple part of an $A_n$-molecule is a weak dual equivalence graph. It remains to show that it 
satisfies the axiom (6), namely that any two $A_{n-1}$-submolecules are connected by a simple edge.

Proceed by induction on $n$, the case $n=1$ being trivial. Let $M$ be an $A_n$-molecule. By inductive assumption, all $A_{n-1}$-molecules 
are Kazhdan-Lusztig. So, according to Lemma \ref{sami:surjective map} there is a covering $M\to \overline{G_\lambda}$, for some partition 
$\lambda$ of $n+1$.

Choose two of these $A_{n-1}$-submolecules of $M$, $A$ and $Z$. Choose a path of simple edges between them which goes through the fewest number of submolecules. 
If it does not go through other submolecules, then we are done. Suppose that is not so. Let $A$, $B$, $C$ be 
the first three submolecules on the path (it may happen that $Z=C$). The partitions $\mu, \nu, \eta$ corresponding to $A,B$, and $C$ are formed by 
removing outer corners of $\lambda$; they are all distinct by Lemma \ref{sami:surjective map}. 

Consider the following string of submolecules connected by simple edges: $A - B - C - A' - B'$, with $A\cong A'$, $B\cong B'$, and 
some of these possibly equalities (this is possible by Lemma \ref{sami:surjective map}). Out of $\mu,\nu$, and $\eta$ choose the partition which is formed by removing the highest box of $\lambda$. 
In the above string, choose a copy of the corresponding submolecule with submolecules attached on both sides (for example, 
if $\lambda\setminus\mu$ was highest of the three, then we should choose $A'$). Then the triple consisting of this submolecule and the two 
adjacent ones satisfies the condition of the Lemma \ref{lemma:shortening} (in the example, it would be the triple $C - A' - B'$). Applying the lemma we get 
that $A' = A$, and $B' = B$. But then $A$ is connected to $C$, contradicting our assumption that the path went through a minimal number of submolecules.

So any two $A_{n-1}$-submolecules are connected by an edge, finishing the proof.
\end{proof}
\end{thm}

\begin{rmk}
In \cite{roberts_degs}, Roberts gives a revised version of Assaf's axiom 6 which is more suitable for computer calculations. Proving our theorem using this alternate axiomatization amounts to 
checking that all the $A_5$-molecules are Kazhdan-Lusztig. This provides a simple computerized proof of our result.
  
\end{rmk}
\newpage
\bibliographystyle{amsalpha}
\bibliography{wgraphs}
\label{sec:biblio}

\end{document}